\documentclass[a4paper,10pt]{amsproc}
\usepackage{amsfonts,lingmacros,tree-dvips, amsmath, amssymb, graphicx, mathrsfs}

\usepackage[all]{xy}

\newdir{ >}{!/6pt/@{ }*:(1,0.2)@_{>}*:(1,-0.2)@^{>}}

\theoremstyle{plain}

\newtheorem{theorem}{Theorem}%[section]

\theoremstyle{definition}
\newtheorem{definition}[theorem]{Definition}

\newtheorem{counter example}[theorem]{Counter Example}

\numberwithin{equation}{section}

\DeclareMathAlphabet{\mathscr}{OT1}{pzc}{m}{it} 
\begin{document}
\Large{
		\title{A FEW REMARKS ON NON-BAIRE SETS IN CATEGORY BASES}
		
		\author[S. Basu]{Sanjib Basu}
		\address{\large{Department of Mathematics,Bethune College,181 Bidhan Sarani}}
		\email{\large{sanjibbasu08@gmail.com}}
		
		\author[A. Deb Ray]{Atasi Deb Ray}
		\address{\large{Department of Pure Mathematics, University of Calcutta, 35, Ballygunge Circular Road, Kolkata 700019, West Bengal, India}}
		\email{\large{debrayatasi@gmail.com}}
		
		\author[A.C.Pramanik]{Abhit Chandra Pramanik}
		\address{\large{Department of Pure Mathematics, University of Calcutta, 35, Ballygunge Circular Road, Kolkata 700019, West Bengal, India}}
		\email{\large{abhit.pramanik@gmail.com}} 
		
		\thanks{The third author thanks the CSIR, New Delhi – 110001, India, for financial support}
	\begin{abstract}
	In this paper, we first establish some equivalent formulations of non-Baire sets in category bases. We then introduce the notion of an uniform non-Baire family of sets and show that there is an uniform non-Baire family inducing a decomposition of the whole space. This phenomenon is then interpreted in the context of the famous Banah-Mazur game.
	\end{abstract}
\thanks{}
	\maketitle

\section*{PRELIMINARIES AND RESULTS}
The idea of category base is a generalization of both measure and topology and its main objective is to present measure and Baire category(topology) and also some other aspects of point set classification within a common framework. It was introduced by J.C.Morgan II in the mid seventies of the last century and has developed since then through a series of papers $[1],[2],[3],[4],[5],$ etc.
To start with, we recall some basic definitions and theorems which may be found in the above references and also in the monograph [6].
\begin{definition}
A category base is a pair (X,$\mathcal{C}$) where X is a non-empty set and $\mathcal{C}$ is a family of subsets of X, called regions satisfying the following set of axioms:
\begin{enumerate}
\item Every point of X belongs to some region; i,e., X=$\cup$$\mathcal{C}$.
\item Let A be a region and $\mathcal{D}$ be a non-empty family of disjont regions having cardinality less than the cardinality of $\mathcal{C}$.\\
i) If A$\cap$($\cup$$\mathcal{D}$) contains a region, then there is a region D$\in$$\mathcal{D}$ such that A$\cap$D contains a region.\\
ii) If A$\cap$($\cup$$\mathcal{D}$) contains no region, then there is a region B$\subseteq$A that is disjoint from every region in $\mathcal{D}$.
\end{enumerate}
\end{definition}

Several examples of category bases can be found in [6].
\begin{definition}
In a category base (X,$\mathcal{C}$), a set is called singular if every region contains a subregion which is disjoint from the set. Any set which can be expressed as countable union of singular sets is called meager. Otherwise, it is called abundant.
\end{definition}

If countable sets are meager, then the category base is point-meager base and a category base whose regions are abundant sets is a Baire base.
\begin{definition}
	In a category base (X,$\mathcal{C}$), a set S is called Baire if in every region, there is a subregion in which either S or its complement X$-$S is meager.
\end{definition}
\begin{theorem}
	The intersection of two regions either contains a region or is a singular sets. 
\end{theorem}
\begin{theorem}(The Fundamental Theorem)
	Every abundant set in a category base is abundant everywhere in some region. This means that for any abundant set A, there exists a region C in every subregion D of which A is abundant.\\
\end{theorem}
In [6], one can find several examples corresponding to Definition 2 and Definition 3.\\

From Definition 3, it follows that a set is non-Baire if there is a region (which is not necessarily unique) in every subregion of which both the set and its complement are abundant. In this paper, we identify this special region as a test region corresponding to the non-Baire set.\\

Below we establish two alternative characterizations for non-Baire sets in a category base. We define 

\begin{definition}
	A set A in a category base (X,$\mathcal{C}$) as having the ($\star$)-property  if there exists a region C such that for every subregion D of C and every set E which essentially contains D , A$\cap$E is non-Baire.\\
	(A set E contains D essentially means that D$-$E is meager [6]) 
\end{definition}
\begin{theorem}
	In any category base, a set is non-Baire if and only if it  satisfies the ($\star$)-property.
	\begin{proof}
		Suppose A is a non-Baire set in a category base (X,$\mathcal{C}$) and C be a test region for A. Let D be any subregion of C. It suffices to show that D can be considered as a test region for A$\cap$P where P is a subset of X which essentially contains D. To prove our contention, we choose any subregion E of D and note that in E both A$\cap$P and its complement are abundant. This is so because in E both E$\cap$A$\cap$P which can be expressed as E$\cap$A$-$A$\cap$(E$-$P) and its complement which is (E$-$A)$\cup$(A$\cap$(E$-$P)) are abundant and E being a subregion of D is also a subregion of C.\\
		
		Conversely, if A is a Baire set then every region C contains a subregion D in which either A or its complement is meager. If D$\cap$A is meage, choose D$-$A as our set E which essentially contains D althouth A$\cap$E=$\phi$ is a Baire set. If D$-$A is meager, choose D$\cap$A as our set E which essentially contains D but A$\cap$E=D$-$A is meager and hence Baire.\\
		
		This proves the theorem. 
	\end{proof}
\end{theorem}

\begin{definition}
	In a category base (X,$\mathcal{C}$), a set A is called completely non-Baire in a region D if for every Baire set B such that B$\cap$D is abundant, both A$\cap$B and (D$-$A)$\cap$B are abundant.
\end{definition}

The above definition is somewhat analogous to the notion of a completely I-nonmeasurable set, given in [8].
\begin{theorem}
	In any Baire base, a set is non-Baire if and only if it is completely non-Baire in some region.
	\begin{proof}
		Let A be a non-Baire set in a Baire base (X,$\mathcal{C}$). Then according to the definition, there is a region D in $\mathcal{C}$ in every subregion of which both A and its complement (X$-$A) are abundant. Let B be any Baire set such that B$\cap$D is abundant. Since by Theorem 4, every region is here a Baire set, we may assume that B$\subseteq$D. From the Fundamental theorem, there is a subregion C of D which is essentially contained in B. Consequently, both B$\cap$A and B$\cap$(D$-$A) are abundant.\\
		
		Conversely, suppose A is completely non-Baire in a region D. Since by hypothesis every region is an abundant Baire set, both A and its complement are abundant in every subregion of D. Hence A is non-Baire.
	\end{proof}
\end{theorem}

Next we introduce 
\begin{definition}
	A family $\mathcal{F}$ of non-Baire sets as an uniform non-Baire family in a category base  (X,$\mathcal{C}$) if there is a region in $\mathcal{C}$ which can act as a common test region for all members of $\mathcal{F}$.
\end{definition}
 With reference to this definition, we prove the following decomposition theorem.
\begin{theorem}
	 If (X,$\mathcal{C}$) is a point-meager, Baire base where $\mathcal{C}$ satisfies CCC (countable chain condition) and every region has cardinality $\aleph_1$, there exists an uniform non-Baire family which induces a decomposition of X.
	 \begin{proof}
	 	Let S be a non-Baire set in (X,$\mathcal{C}$) and T=X$-$S. The existence of S is ensured by Theorem 6, Ch 2,Sec II, [6]. Then there is a region D in every subregion of which both S and T are abundant. Now by Grzegorek's unification of Sierpinski's theorems [1], (see also Theorem 16, Ch 2, Sec II, [6]) both S and T can be decomposed into families \{$S_\alpha\}_{\alpha < \Omega}$ and \{$T_\alpha\}_{\alpha < \Omega}$ ($\Omega$ is the first ordinal of cardinality $\aleph_{1}$) of mutually disjoint sets such that each $S_\alpha$ (resp. $T_\alpha$) is abundant in every region in which S (resp. T) is abundant. Since $X-S\subseteq X-S_\alpha$ (resp. X$-T\subseteq X-T_\alpha$), so every S$_\alpha$ and its complement (resp. T$_\alpha$ and its complement ) are abundant in every subregion of D proving that S$_\alpha$ (resp. T$_\alpha$) are both non-Baire sets for all $\alpha < \Omega$ and D can act as a common test region for the entire family \{$S_\alpha$, $T_\alpha$ : $\alpha < \Omega$\}. Moreover, this family constitutes a partition of X. Thus we find the existence of an uniform non-Baire family inducing a decomposition of X.
	 	\end{proof}
\end{theorem}
A generalization of the classical Banach-Mazur game [7] was presented in [3] by Morgan through the introduction of the notion of $\mathcal{M}$-family of sets, which was defined in the following manner : \\
\begin{definition}
	 A family $\mathcal{C}$ of subsets of a non-empty set X is called an $\mathcal{M}$-family if it satisfies the following set of axioms :
	\begin{enumerate}
		\item The intersection of any descending sequence of $\mathcal{C}$-sets is non-empty.
		\item Suppose x is a point in X, then\\
		i) there is a $\mathcal{C}$-set  containing x, i.e., X=$\cup$$\mathcal{C}$.\\
		ii) for each $\mathcal{C}$-set A, there is a $\mathcal{C}$-set B$\subseteq$A such that x$\notin$B.
		\item Let A be a $\mathcal{C}$-set and D be a nonempty family of disjoint $\mathcal{C}$-sets having cardinality less than the cardinality of $\mathcal{C}$.\\
		i) If A$\cap$($\cup$$\mathcal{D}$) contains a $\mathcal{C}$-set , then there is a $\mathcal{D}$-set D such that A$\cap$D contains a $\mathcal{C}$-set.\\
		ii) If A$\cap$($\cup$$\mathcal{D}$) contains no $\mathcal{C}$-set, then there is a $\mathcal{C}$-set B$\subseteq$A which is disjoint from every set in $\mathcal{D}$.
	\end{enumerate} 
\end{definition}
Evidently, any $\mathcal{M}$-family is a special case of a category base, so the definitions of meager, abundant and Baire sets are as usual. Moreover, any $\mathcal{M}$-family is point-meager and every member of it is an abundant set which follows directly from the axioms (1) and (2) in the above Definition.\\

If $\mathcal{C}$ is a non-empty family of subsets of a non-empty set X and S is a subset of X, then as described in [3], the generalized Banach-Mazur game $\Gamma$(S, $\mathcal{C}$) is played as follows : Two players I and II, alternatively choose sets from $\mathcal{C}$ to define a descending sequence of sets, player I selecting the sets in the sequence with odd indices and player II selecting the sets with even indices. If the intersection of the constructed sequence has atleast one point (which obviously exists if $\mathcal{C}$ is an $\mathcal{M}$-family) in S, then player I wins; otherwise player II wins. In the game $\Gamma$(S, $\mathcal{C}$), let us denote player I, player II by using symbols $\langle$S$\rangle$, $\langle$X$-$S$\rangle$ respectively.\\

Now suppose C$\in$$\mathcal{C}$ such that S$\cap$C $\neq$ $\phi$ $\neq$ (X$-$S)$\cap$C. Then we can describe the game $\Gamma$(S,$\mathcal{C}_{|C}$) in a manner similar as above where $\mathcal{C}_{|C}$=\{E$\in$$\mathcal{C}$ : E$\subseteq$C\}. As a consequence of the Theorem 2 [3], Theorem 11 may be interpreted in terms of the generalized Banach-Mazur game in the following manner :
\begin{theorem}
	 Let $\mathcal{C}$ be an $\mathcal{M}$-family satisfying CCC (countable chain condition) and every $\mathcal{C}$-set has cardinality $\aleph_{1}$. Then there exist a $\mathcal{C}$-set D and a partition X=$\bigcup\limits_{\alpha < \Omega} {X_\alpha}$  such that under the condition that there exists a sequence $\{h_n\}_{n=1}^{\infty} (h_n :\mathcal{C}\mapsto\mathcal{C})$ satisfying
	 \begin{enumerate}
	 	\item for every $\mathcal{C}$-set A, h$_n$(A)$\subseteq$A ;
	 	\item for every sequence $\{A_n\}_{n=1}^{\infty}$ of   $\mathcal{C}$-sets, if $\{h_n(A)\}_{n=1}^{\infty}$ is descending, then $\bigcap\limits_{n=1}^{\infty}h_n$(A) contains only one point;
	 \end{enumerate}
	  no player  $\langle X_\alpha \rangle$ ($\alpha < \Omega$) can have a winning strategy in the game $\Gamma(X_\alpha,\mathcal{C}_{|D}$).
\end{theorem}

%\textbf{Acknowledgments: } The authors would like to express their gratitude to the learned referees who have given numerous valuable suggestions towards the improvement of the initial version of this article. 
\bibliographystyle{plain}

	\end{document}